\documentclass[%
 reprint,
groupedaddress,
nofootinbib,
 amsmath,amssymb,
 aps,
 pre,
]{revtex4-2}

\usepackage{dcolumn}
\usepackage{bm}

\usepackage{amsmath}
\usepackage{amssymb}
\usepackage{amsthm}
\usepackage{dsfont} 
\usepackage{float}
\usepackage{multirow}
\usepackage{siunitx}
\usepackage{amssymb}
\usepackage{amsmath}
\usepackage{mathrsfs}
\usepackage{mathtools}
\usepackage{stmaryrd}
\usepackage{siunitx}
\usepackage{booktabs}
\usepackage{tabularx}
\usepackage{bbm}
\usepackage{adjustbox}
\usepackage{threeparttable}
\usepackage[english]{babel}
\usepackage{algpseudocode}
\usepackage{algorithm}

\newtheorem{thm}{Theorem}
\newtheorem{defi}{Definition}
\newtheorem{prop}{Proposition}
\newtheorem{lem}{Lemma}
\newtheorem{rmk}{Remark}
\newtheorem{nota}{Notation}
\DeclareMathOperator*{\argmax}{argmax}

\usepackage{mathtools, stmaryrd}
\usepackage{xparse} \DeclarePairedDelimiterX{\Iintv}[1]{\llbracket}{\rrbracket}{\iintvargs{#1}}
\NewDocumentCommand{\iintvargs}{>{\SplitArgument{1}{,}}m}
{\iintvargsaux#1} %
\NewDocumentCommand{\iintvargsaux}{mm} {#1\mkern1.5mu..\mkern1.5mu#2}

\newcommand{\innerproduct}[2]{\langle #1, #2 \rangle}

\begin{document}

\preprint{APS/123-QED}

\title{Generating social networks with static and dynamic utility-maximization approaches}
\thanks{We thank F. Tarissan and F. Bloch for fruitful and useful discussions and J. Randon-Furling for his advice.}

\author{Aldric Labarthe}
 \email{aldric.labarthe@ens-paris-saclay.fr (corresponding author)}
\author{Yann Kerzreho}%
 \email{yann.kerzreho@ens-paris-saclay.fr}
\affiliation{
 Ecole Normale Supérieure (ENS) Paris-Saclay, 
Gif-sur-Yvette, 91190 France
}

\date{\today}

\begin{abstract}
In this paper, we introduce a conceptual framework that model human social networks as an undirected dot-product graph of independent individuals. Their relationships are only determined by a cost-benefit analysis, \textit{i.e.} by maximizing an objective function at the scale of the individual or of the whole network. On this framework, we build a new artificial network generator in two versions. The first fits within the tradition of artificial network generators by being able to generate similar networks from empirical data. The second relaxes the computational efficiency constraint and implements the same micro-based decision algorithm, but in agent-based simulations with time and fully independent agents. This latter version enables social scientists to perform an in-depth analysis of the consequences of behavioral constraints affecting individuals on the network they form. This point is illustrated by a case study of imperfect information. 

\end{abstract}
\maketitle

\section{\label{sec:intro} Introduction}

Existing network models often fall short in accurately representing real social networks, primarily due to limitations in providing a consistent and flexible modeling framework. By consistent we mean the ability to have, at every scale (node or whole network), rules that define agents interactions in a unified way without any \textit{ad hoc} mathematical or computational processes designed to meet pre-established requirements, \textit{e.g.} by (re)-defining edges randomly to check the small world property or to satisfy statistical requirements \cite{watts1998collective, barabasi1999emergence, jackson2007meeting, hamill2009social, pham2012s3g2, erling2015SNB, zhuge2018agent}. This excludes the latest advances in graph reproduction using neural networks \cite{liao2019efficient, zhu2022survey}, which are black boxes with no consistent rules on graph formation. By flexible we mean a framework that can cover the whole set of possible human networks and that can be easily modified and extended, \textit{e.g.} towards agent-based modeling.

In the field of Random Geometric Graphs (RGG), nodes are points in a latent space, and their distances from each other allow to decide independently whether they are connected or not (for a review see \cite{duchemin2023random}). The geometry gives good clustering properties for modeling social networks \cite{krioukov2016clustering} but as a distance must verify the triangular inequality, the small-world property is difficult to achieve \cite{talaga2019homophily, iijima2017, papamichalis_latent_2022, saha_study_2024}. The RGG framework is also flexible as different link spaces and functions can be used, and implementing an agent-based model approach is possible \cite{gonzalez2006networks, caux2014dynamic} with real data \cite{klopp2017oracle, klopp2019optimal}. However, the independence assumption between edges is an important limitation because it does not take into account the presence of possible substitution effects between connections. In this paper, we introduce a different model that could remove the independence assumption and avoid the geometric constraints induced by distance. We indeed place ourself in continuation with dot product graphs models \cite{young_directed_2008} by using the dot product of the latent space instead of the distance to assess a proximity or compatibility score between individuals. 

Following the seminal work of \cite{jackson1996} several approaches in the economics literature have tried to model social networks with agent-based postulates, where agents maximize their utility by managing their edges which generate costs and benefits (for an extended review, we refer to \cite{hellmann2014, jackson2007study, jackson2011overview}). This part of the literature uses linear relationship costs and benefits that depend on the graph and not on individuals. \cite{hellmann2014} stated that the benefits of a relationship are always assumed to be geometrically linked to a closeness measure (\cite{jackson1996}, or \cite{bonacich1987power} for another definition of closeness). This lack of true link between relationship benefits and specific individuals characteristics leads to contradictions, \textit{e.g.} no individual has a rational incentive in becoming a hub and having many more links than other individuals. Additionally, technical considerations are still quite overlooked in the economics literature \cite{jackson2007study}. Almost all models heavily rely on agent-based simulations \cite{tesfatsion1997trade} that cannot match the convenience and the efficiency of simulators encountered in the RGG or the machine learning literature. 

Our approach can be seen as micro-based, as we define a maximization problem in which each weighted-undirected edge (relationship between individuals) induces a benefit and a cost. Each pair of individuals shares a compatibility which drives up the benefit from the relationship. The weight of an edge represents the intensity of the relationship, accounting for the diversity of social interactions. Our postulates make possible to create graphs by maximizing an objective function. Depending on the scale of the optimization (individual or social), the objective function is the utility function of the individual, or the sum of the utilities of each individual. Therefore, we introduce two artificial network generators: the first uses an analytical optimization process of the social well-being of the population, the latter takes advantage of a more flexible agent-based framework to generate networks that are the result of individual decisions. Hence, our objectives are twofold: using our first model we will introduce a generator that is able to reproduce, with high fidelity, every empirical networks (providing they meet some constraints). Secondly, we will use the flexibility of our second generator to study the effects on the network of imperfect information by limiting the agent's scope of action.

We begin by introducing the mathematical framework on which rely both of our generators, and theoretically define our conceptual framework and generators in Section \ref{sec:methodology}. The results are then presented in Section \ref{sec:results} and discussed in Section \ref{sec:discussion}.

\section{Methodology} \label{sec:methodology}

\subsection{Conceptual Framework}

\begin{defi}[Individuals space]\label{def_indivSpace}
We denote as \textbf{Individuals space} the product space $$\mathcal{I} = (\mathcal{E}, \innerproduct{.}{.}) \times \mathcal{U}_I$$ with $(\mathcal{E}, \innerproduct{.}{.})$ a finite-dimensional euclidean space defined on the vector space $\mathcal{E}$ with $\innerproduct{.}{.}$ an inner product over $\mathbb{R}^{+}$. Individuals sets are elements of $\mathcal{I}^I$. No individual should be orthogonal to any other in its first component.
\end{defi}

\begin{defi}[Individuals graph]\label{def_indivGraph}
For a given individuals space $\mathcal{I}$ and a given individuals set $\mathscr{I} \in \mathcal{I}^I$, an \textbf{individual graph} is a 3-tuple $\mathcal{G}(\mathscr{I} , \alpha, C)$ with $\alpha, C$ two vectors of $(\mathbb{R}^{*}_{+})^{n_I}$ and $(\mathbb{R}_{+})^{n_I}$) with $n_I=\binom{I}{2} = \frac{I^2-I}{2}$ :
\begin{itemize}
    \item $\alpha$ is a \textbf{weight list}. $$\alpha = (\alpha_{1,2}, \hdots, \alpha_{1,I},  \alpha_{2,3}, \hdots, \alpha_{2,I}, \hdots,\alpha_{I-1,I} )^T$$
    with $\alpha_{i,j}$ the edge weight between agent i and j
    \item $\chi(i,j) = \innerproduct{P_i}{P_j}\mathbbm{1}(i\neq j)$ is the \textbf{compatibility function}, with $(i,j)\in \mathscr{I}^2$ and $P_i$ the first component of i (in $\mathcal{E}$).
    \item $C$ is a \textbf{compatibility list} $$C = (c_{1,2}, \hdots, c_{1,I},  c_{2,3}, \hdots, c_{2,I}, \hdots,c_{I-1,I} )^T$$ with $c_{i,j} = \chi(i,j) = \innerproduct{P_i}{P_j}$.
    
\end{itemize}

On each graph is defined $d : \mathcal{I}^2 \rightarrow \mathbb{N}$ the distance between agents $i$ and $j$: the smallest number of intermediary nonnull edges it takes to reach a vertex from another.
\end{defi}

\begin{nota}\label{nota_cAndAlpha}
We denote $\forall i \in \Iintv*{1, I }, \alpha^i = (((\alpha_{i,j})_{i<j\leq I}), ((\alpha_{j,i})_{1\leq j<i}))$ the (I-1)-dimensional vector of weights associated with all possible edges linking the vertex $i$ with other individuals. %$\alpha_i$ is the ith-component of $\alpha$, and 
$\alpha^i_j$ is the jth-component of $\alpha^i$. When we write $\alpha_{i,j}$, we consider $\alpha_{\min(i,j), \max(i,j)}$. We use the same notations for $C$ (with $C^i$ and $c^i_j$).
\end{nota}

If we implement a latent space as usual in the RGG field \cite{iijima2017}, we do not use a standard distance between individuals but rather the dot product of their characteristic vectors. This choice gives us more freedom, as a standard metric distance would impose compatibilities to satisfy the triangular inequality \cite{iijima2017,talaga2019homophily} and would prevent us from finding a suitable compatibility list for every empirical network. In addition, we choose not to provide a probability framework as commonly done in the RGG literature as the focus of this study is essentially placed on empirical networks, and not theoretical ones. However, one can observe that it is perfectly possible to consider a prior distribution for the $P$ variables. 

\begin{defi}[ALKY utility function]\label{def_ALKYUtility}
We consider the set $\mathcal{U}^{\text{ALKY}}_I$ of applications $U : [0,1)^{I-1} \times (\mathbb{R}^{*}_{+})^{I-1} \longrightarrow \mathbb{R}$:
\begin{equation*}
    \begin{split}
\mathcal{U}^{\text{ALKY}}_I = \left\{\kappa>1, \gamma>1, \delta>1, \right. \\ U\left(\alpha, c \right)= \left. \sum^{I-1}_{j=1} \left[ \kappa \alpha_{j} c_{j} - \frac{\alpha_{j}^{\gamma}}{1-\alpha_{j}^{\gamma}} \right] - \left( \sum^{I-1}_{j=1} \alpha_{j} \right)^{\delta}\right\}    
\end{split}
\end{equation*}

\end{defi}

The ALKY function is made up of three parts: a positive one which weights  the affinity $(c_{i, j})_{j \leq N}$ by the relationship intensity $(\alpha_{i, j})_{j \leq N}$ and two costs, that convexly penalize intensity, and create substitution effects between relationships. $\gamma$ is the convexity parameter of relationships' direct cost (and prevents corner solutions). $\delta$ enable the substitution effects between relationships. For $\delta > 1$ relations are substitutable: an agent prefers reducing the emphasis on relationships with lower affinity in order to increase the weight on his most successful matches. 

\begin{rmk}
All results are stated for a more generic class of utility functions that we define in the supplementary materials.
\end{rmk}

\subsection{Optimal Graph generation: a social optimization}

\begin{defi}[Social Optimization Problem]\label{def_SOP}
For a given individuals graph $\mathcal{G}_N(\mathscr{I} , \alpha, C)$, we define the Social Optimization Problem (SOP) as:
\begin{equation*}
\begin{cases}
\argmax_{\alpha} \mathcal{O} (\alpha, C) =  \sum^N_{i=1} U_i(\alpha^i, C^i)\\
\text{s.t} \quad \begin{cases} 
\alpha \in [0,1)^{n_N} \\
\nexists i \in \Iintv*{1, N }, \forall j \in \Iintv*{1, N-1 }, \alpha^i_j = 0
\end{cases}
\end{cases}
\end{equation*}
{\footnotesize With $U_i$ the utility function of agent $i$, that is its second component in the individuals space.}
\end{defi}

\begin{thm}[Unicity of the optimal graph]\label{thm_unicity}
For any $\mathcal{G}(\mathscr{I} , \alpha, C)$, the Social Optimization Problem has only one unique solution.
\end{thm}

We denote as \textit{Optimal Graph} the unique solutions of the Social Optimization Problem introduced in Definition \ref{def_SOP}. Using Theorem \ref{thm_alphaIsAlwaysOpti}, we know that every network that is associated with a weighted and undirected graph, without any isolated vertex nor loop, can be generated using our model. In addition, all generated networks verify pairwise stability \cite{jackson1996} (see the supplementary materials for a proof).

\begin{thm}\label{thm_alphaIsAlwaysOpti}
For any individuals set $\mathscr{I} \in (\mathcal{E}, \innerproduct{.}{.}) \times \mathcal{U}_N$, and weight list $\alpha$ of dimension $n_N = \binom{N}{2}$ satisfying:
\begin{itemize}
\item $\forall a \in \alpha, a \in [0,1)$
\item $\nexists i \in \Iintv*{1, N }, \forall j \in \Iintv*{1, N-1 }, \alpha^i_j = 0$
\item $\forall i \in \Iintv*{1, N }, (\mathscr{I})_i = (P_i, U_i), U_i \in \mathcal{U}^{\text{ALKY}}_I\;$ \footnote{This condition can be loosen by using the more general definition of utility functions, see the supplementary materials for more details.}
\end{itemize}

There exists:
\begin{enumerate}
\item a unique compatibility list $\bar C$ such that $\alpha$ is the unique unconstrained solution of the Social Optimization Problem for the individuals graph $\mathcal{G}(\mathscr{I} , \alpha, \bar C)$
\item a set of compatibility lists $\underbar C$ such that $\forall C \in \underbar C$, $\alpha$ is the unique solution of the Social Optimization Problem for the individuals graph $\mathcal{G}(\mathscr{I} , \alpha, C)$. 
$$\underbar C = \left\{ C \in (\mathbb{R}^{*}_{+})^{n_I}, \quad c_{ij} \in \begin{cases} \{\bar c_{ij}\} & \text{if } \alpha_{ij} > 0\\  (0, \bar c_{ij}]& \text{if } \alpha_{ij} = 0\end{cases}\right\}$$
\end{enumerate}
\end{thm} 

\begin{figure}[h]
    \centering
    \includegraphics[width=0.45\textwidth]{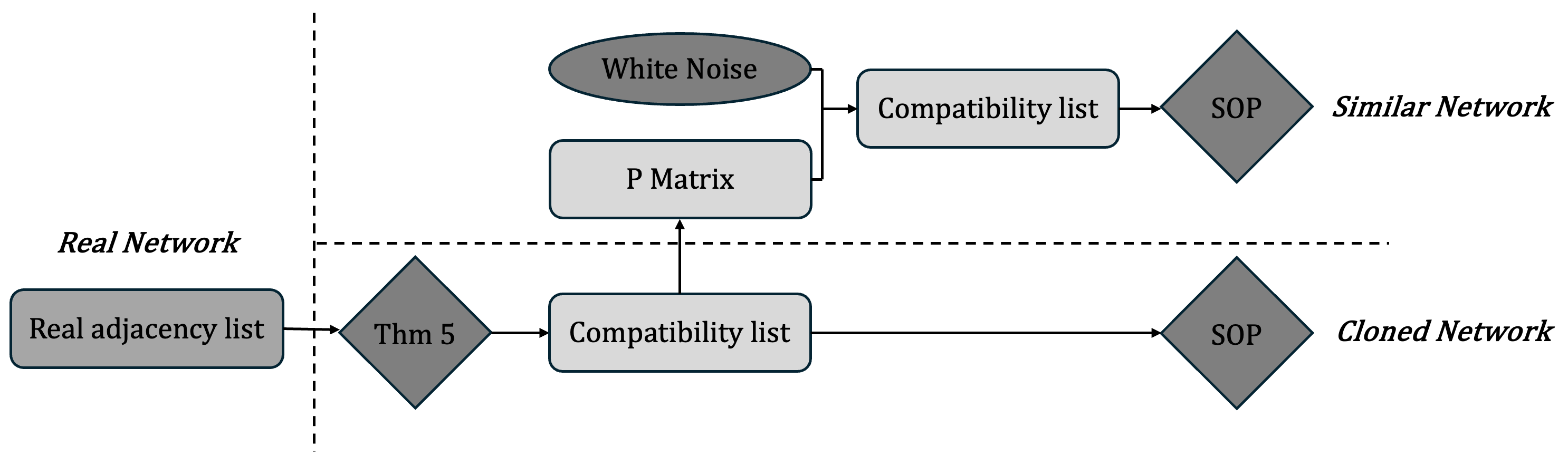}
    \caption{\textbf{Methodology to generate similar networks or perfect clones of a real network} }  \label{fig:clonage}
\end{figure}

Moreover, this property will allow us to derive multiple similar synthetic graphs from any real graph, that will share some overall benchmark features, but with different edges and weights. Our methodology is illustrated in Figure \ref{fig:clonage}: we first use Theorem \ref{thm_alphaIsAlwaysOpti} to recover the compatibility matrix associated with a given tuple of parameters (that can be freely chosen) from the real adjacency list. With this tuple of parameters and this compatibility list, the model is always able to regenerate cloned networks (as the solution of the SOP is unique). 

We name P Matrix in figure \ref{fig:clonage} the matrix that is formed by all individuals' characteristic vectors: the ith row of the P Matrix is the characteristic vector of the individual i (in $\mathcal{E}$). Several linear algebra methods can be used to recover candidates for these characteristic vectors. A proof of existence of such a P Matrix is provided in the supplementary materials, and setting the dimension of the latent space to the number of individuals can be taken as a rule of thumb. We can alter this matrix with white-noise processes with different standard deviations. This leads to a new compatibility list and, by applying the SOP to this list, we finally obtain a new graph derived from the real graph.

Theoretically, the resulting network from the Social Optimization problem (SOP) can be seen as a social equilibrium (as usual in the economics literature), that is what would be chosen to maximize the global welfare of the network of individuals. In this approach, all individuals are perfectly substitutable, i. e. it can be optimal to reduce the welfare of one individual to generate a higher welfare for another. One can argue that the model lacks some real properties. 

First, the substitution principle does not seem to be appropriate as some individuals can refuse a relationship that does not bring them strictly positive utility gains, and therefore refuse to reduce their welfare for the greater good. Second, both social (Definition \ref{def_SOP}) and individual (Definition \ref{def_IOP}) optimization problems are instantaneous and timeless. Third, in both settings, information (in our model, the knowledge of every individual in the network) is assumed to be freely and fully available to all individuals, which is, in the context of human networks, a costly hypothesis, considering that an individual cannot know its affinity with every single individual available. Hence, we will now study networks that are not generated using our Social Optimization Problem (Definition \ref{def_SOP}) but that are the result of a decentralized optimization problem.

\subsection{Decentralized graph generation: an agent-based model}
To tackle these two limitations, we introduce an agent-based model that tries to implement the Individual Optimization Problem (Definition \ref{def_IOP}) with temporal dynamics and limited information. 

\begin{defi}[Individual Optimization Problem]\label{def_IOP}
For a given individuals graph $\mathcal{G}_N(\mathscr{I} , \alpha, C)$, we define for all individual $i \in \Iintv*{1, N }$ the Individual Optimization Problem (IOP) as:
\begin{equation*}
\begin{cases}
\argmax_{\alpha^i} U_i(\alpha^i, C^i) \\
\text{s.t.} \quad \alpha^i \in [0,1)^{I-1}
\end{cases}
\end{equation*}
\end{defi}

\begin{thm}[Unicity of the individual optimal graph]\label{thm_unicityIOP}
For any $\mathcal{G}(\mathscr{I} , \alpha, C)$, and for all individual $i \in \Iintv*{1, N }$ the Individual Optimization Problem (IOP) has only one unique solution.
\end{thm}

Algorithms \ref{alg:agent} and \ref{alg:agent_response} provide details about the implementation: at each time step, each individual is asked (in a random order) to select an action according to the gradient of its utility function. This action can be either to create a new relationship (hereinafter, edge), to increase or decrease the weight of an existing edge, or to delete a edge. We introduce the concept of scope, which is a set containing all individuals that a given agent is able to interact with. The scope of an individual $i$ is defined as $\{j : d(j, i) \leq \delta \}$, with $\delta$ the scope depth. If the agent wants to increase or create an edge, he has to ask the other individual for his approval, which is again evaluated using the gradient of its own utility evaluated on his own scope. This approach guarantees reciprocity, and implements both constraints of time and imperfect information, as agents are evaluating the gradient of their utility function not on the whole network but only on their scope. Moreover, as an exploration process, agents are selecting their actions using softargmax. A memory of the 3 last actions (and requests) is implemented, preventing individuals to ask again the same others for 3 consecutive periods, and thus incentivizing individuals to explore and select actions that have not the highest gradient value, as it does not take into account possible future gains induced by the scope increase.

\begin{figure}[h]
\begin{algorithm}[H]
\caption{Individual $i$ select action protocol (simplified)}\label{alg:agent}
\begin{algorithmic}

\State $\mathcal{S}^i \gets \{j : d(j, i) \leq \delta \}$ \Comment{Get individuals in the scope of $i$}

\State $C^i_\mathcal{S} \gets \{c_{i,j} : j \in \mathcal{S}^i \}$ \Comment{Get the set of compatibilities}

\State $\alpha^i_\mathcal{S} \gets \{\alpha_{i,j} : j \in \mathcal{S}^i \}$ \Comment{Get the set of weights}

\State $k \gets \text{softargmax}_j \left|\nabla U^i(\alpha^i_\mathcal{S}, C^i_\mathcal{S})\right|$ 

\If{$\nabla U^i_{(k)} > 0$}
    \State \Call{request}{$k, i, \alpha_{i,k} + \lambda \nabla U^i_{(k)}$} 
\Else
    \State $\alpha_{i,k} \gets \alpha_{i,k} + \lambda \nabla U^i_{(k)}$ \Comment{Direct decrease in friendship with $N^k$}
\EndIf

\end{algorithmic}
\end{algorithm}
\begin{algorithm}[H]
\caption{Individual $k$ response protocol to $i$ (request function, simplified)}\label{alg:agent_response}
\begin{algorithmic}
\Function{request}{$k, i, a$}
    \State $\mathcal{S}^k \gets \{j : d(j, k) \leq \delta \} \bigcup \{i\}$ \Comment{Get individuals in scope}
    
    \State $C^k_\mathcal{S} \gets \{c_{k,j} : j \in \mathcal{S}^k$ \} \Comment{Get the set of compatibilities}
    
    \State $\alpha^i_\mathcal{S} \gets \{\alpha_{k,j} : j \in \mathcal{S}^k$ \} \Comment{Get the set of weights}
        
    \If{$\nabla U^k_{(i)}(\alpha^k_\mathcal{S}, C^k_\mathcal{S}) > 0$}
        \State $\alpha_{i,k} \gets \text{min}(a, \alpha_{i,k} + \lambda \nabla U^k_{(i)})$
        \State \Return true
    \EndIf
    \State \Return false
\EndFunction

\end{algorithmic}
\end{algorithm}
\end{figure}

Theorem \ref{thm_invisibleHand} provides some conditions on the compatibility structure for the equilibrium to exists, considering that the equilibrium exists if and only if for each pair of individuals $(i,j)$, they both desire the same $\alpha_{i,j}$ at equilibrium. Under these conditions, the theoretical individual optimal graph are the same as the social optimization problem. This setting will allow us to compute the theoretical optimal graph with our standard generator, and then compare the networks obtained with our decentralized generator (ABM) with respect to it. 

\begin{thm}[Invisible hand]\label{thm_invisibleHand}
For any $\mathcal{G}(\mathscr{I} , \alpha, C)$ with $\mathcal{C}$ the compatibility matrix defined as $\mathcal{C} = (c_{i,j})_{1\leq i,j\leq N}$ ($\mathcal{C}^{i}$ its ith-column), and $(\alpha^i)^*$ the solution vector of the Individual optimization problem (IOP) evaluated on $c^i$, the compatibility list derived from $\mathcal{C}_{i}$, if
\begin{enumerate}
\item $\exists \bar c \in \mathbb{R}^{N} \text{s.t. } \forall i \in \Iintv*{1, N }, \mathcal{C}^i \in \mathscr{S}(\bar c)$, with $\mathscr{S}(\bar c)$ the symmetric group of $\bar c$
\item there exists an ALKY utility function such that $\forall i \in \Iintv*{1, N }, (\mathscr{I})_i = (P_i, U_i), U_i = U \;$\footnote{This condition can be loosen by using the more general definition of utility functions, see the supplementary materials for more details.}
\end{enumerate}
Then there exists $\hat \alpha  \in \mathbb{R}^{n_N}$ such that $\forall i \in \Iintv*{1, N }, (\alpha^i)^* \in \mathscr{S}(\hat \alpha)$. Moreover, $\alpha^{**}$ the solution of the associated SOP, verifies: $$\forall (i,j) \in \Iintv*{1, N }^2, \alpha^{**}_{i,j} = (\alpha^{\max{\{i,j\}}})^*_{\min{\{i,j\}}}$$
 \end{thm}

\subsection{Simulation methodology}
Hereinafter, we only use homogeneous individuals sets, where all individuals share the same ALKY utility function from Definition \ref{def_ALKYUtility}, with parameters $\gamma = 9, \delta = 2, \kappa = 10$. This should not have any effect on our results as there always exists a compatibility list that produces the same output network with another choice of parameters, as the model has two degree of freedom. 

For the Decentralized Graphs, we use Theorem \ref{thm_invisibleHand} to create compatibility structures in which the Social optimization problem (Definition \ref{def_SOP}) is equivalent to the decentralized problem, that is the union of Individuals optimization problem (Definition \ref{def_IOP}). We generate 50 random compatibility structures of 64 individuals that meet the requirements of Theorem \ref{thm_invisibleHand}. Although it is perfectly impossible to cover the whole set of all possible compatibility structures, our sample has been designed to contain heterogeneous networks (the global clustering coefficient is distributed between 0.107 and 0.918, the Gini coefficient between 0.003 and 0.309, the density between 0.127 and 0.921). We use stochastic initialization adjacency lists: $\alpha_{i,j} = (\frac{\omega}{2}, 0 ; \frac{\iota}{I}, 1 - \frac{\iota}{I})$, with $\omega$ the minimum non-negative weight allowed and $\iota$ the initialization density. For each compatibility structure, and for each parametrization ($\iota$ and $\delta$), we will perform 3 simulations and provide means and standard deviations. All simulations are performed with a fixed duration, but all results are extracted from simulations that can be considered as having reached a limit point: all indicators (global clustering coefficient, density, Gini, ...) do not exhibit any statistically significant variations in the last 5\% time, or the round-by-round simulation has registered no accepted actions in the last 50 rounds.

\section{Results} \label{sec:results}

\subsection{Similar networks generation}
We here demonstrate on an example the ability of the model not only to regenerate the network, but to enable the production of very similar networks in their overall benchmark properties, but with different edges and weights. We use as an example the real network named "Moreno Train bombing" \cite{konect2017morenotrain} available in the Konect database \cite{konect}. This network is made of "contacts between suspected terrorists involved in the train bombing of Madrid on March 11, 2004 as reconstructed from newspapers."

\begin{table}[] 
\begin{adjustbox}{max width=0.45\textwidth}
  \begin{threeparttable}[t]
    \centering
    \caption{\textbf{Similar networks generation methodologies assessment on the example of the Train bombing network from Konect} (64 individuals). \textit{All values are means (with a precision of 1e$^{-2}$) on a sample of 100 networks generated using the Social Optimization Problem (SOP). Standard deviations across the 100 compatibility structures are in parenthesis with a precision of 1e$^{-3}$.}} \label{tab:compatibilityClonage}
\begin{tabular}{lcccccc}
\midrule\midrule
% & & \multicolumn{7}{c} {Clones} \\
%   \cmidrule(lr){3-9}
 & & \multicolumn{2}{c} {WN($\sigma$) on C} & \multicolumn{3}{c} {WN($\sigma$) on P} \\ 
  \cmidrule(lr){3-4}  \cmidrule(lr){5-7}
 & Original & 0.1 & 0.0001 & 0.2 & 0.1 & 0.05  \\
\midrule
Vertices & 486 & 325.14 & 457.76  & 765.35 & 489.12 & 371.59 \\
 & {\scriptsize(-)} & {\scriptsize(6.571)} & {\scriptsize(5.842)}  & {\scriptsize(65.848)} & {\scriptsize(30.114)} & {\scriptsize(12.861)} \\
Del. Vertices & - & 407.46 & 319.96 & 368.62 & 392.44 & 399.36 \\
 & {\scriptsize(-)} & {\scriptsize(8.394)} & {\scriptsize(9.268)} & {\scriptsize(15.019)} & {\scriptsize(10.28)} & {\scriptsize(9.304)} \\
New Vertices & - & 246.6 & 291.72 - & 647.98 & 395.56 & 284.94  \\
 & {\scriptsize(-)} & {\scriptsize(10.657)} & {\scriptsize(10.483)} & -{\scriptsize(61.387)} & {\scriptsize(29.745)} & {\scriptsize(15.605)}  \\
 \midrule
Clustering & 0.56 & 0.1 & 0.24  & 0.75 & 0.55 & 0.29  \\
 & {\scriptsize(-)} & {\scriptsize(0.016)} & {\scriptsize(0.012)} & {\scriptsize(0.035)} & {\scriptsize(0.06)} & {\scriptsize(0.049)} \\
Gini coef. & 0.13 & 0.19 & 0.23 & 0.12 & 0.18 & 0.21  \\
 & {\scriptsize(-)} & {\scriptsize(0.013)} & {\scriptsize(0.008)} & {\scriptsize(0.013)} & {\scriptsize(0.012)} & {\scriptsize(0.012)}  \\
Density & 0.12 & 0.08 & 0.11 & 0.19 & 0.12 & 0.09 \\
 & {\scriptsize(-)} & {\scriptsize(0.002)} & {\scriptsize(0.001)} & {\scriptsize(0.016)} & {\scriptsize(0.007)} & {\scriptsize(0.003)}  \\
Avg. distance & 0.63 & 1.25 & 0.5  & 1.99 & 1.66 & 1.39 \\
 & {\scriptsize(-)} & {\scriptsize(0.065)} & {\scriptsize(0.02)} & {\scriptsize(0.348)} & {\scriptsize(0.182)} & {\scriptsize(0.122)} \\
 \midrule
Assortativity & 0.03 & -0.05 & -0.11  & 0.15 & 0.09 & -0.18 \\
 & {\scriptsize(-)} & {\scriptsize(0.057)} & {\scriptsize(0.035)} & {\scriptsize(0.13)} & {\scriptsize(0.12)} & {\scriptsize(0.086)}  \\
Avg. degree & 7.59 & 5.08 & 7.15 & 11.96 & 7.64 & 5.81  \\
 & {\scriptsize(-)} & {\scriptsize(0.103)} & {\scriptsize(0.091)} & {\scriptsize(1.029)} & {\scriptsize(0.471)} & {\scriptsize(0.201)}  \\
Var. degree & 37.99 & 7.47 & 30.19 & 43.19 & 16.2 & 10.67  \\
 & {\scriptsize(-)} & {\scriptsize(0.508)} & {\scriptsize(0.876)} & {\scriptsize(6.537)} & {\scriptsize(2.699)} & {\scriptsize(0.819)}  \\
 \midrule
S.d. Comp. & 0.28 & 0.3 & 0.28 & 2.59 & 0.74 & 0.35  \\
 & {\scriptsize(-)} & {\scriptsize(0.002)} & {\scriptsize(0)} & 
 {\scriptsize(0.433)} & {\scriptsize(0.111)} & {\scriptsize(0.02)} \\
Cor.\tnote{1}  Comp. & 0.9 & 0.73 & 0.9  & 0.05 & 0.21 & 0.55  \\
 & {\scriptsize(-)} & {\scriptsize(0.005)} & {\scriptsize(0)} & {\scriptsize(0.02)} & {\scriptsize(0.042)} & {\scriptsize(0.055)}  \\
S.d. Cor.\tnote{2} & 0.14 & 0.12 & 0.14  & 0.89 & 0.69 & 0.36 \\
 & {\scriptsize(-)} & {\scriptsize(0.003)} & {\scriptsize(0)} & {\scriptsize(0.027)} & {\scriptsize(0.043)} & {\scriptsize(0.051)} \\
\midrule\midrule\end{tabular}
     \begin{tablenotes}
       \item [1] Cor. denotes the Pearson coefficient between individuals' compatibility lists (for agent i, $c^i=(c_{1,i}, ..., c_{N,i})$)
       \item [2] Sd. Cor. is the standard dispersion of the Pearson coefficients.
     \end{tablenotes}
  \end{threeparttable}
  \end{adjustbox}
\end{table}

One naive approach would be to try to directly alter the compatibility list. Nonetheless the model's response to compatibility list is not continuous. Table \ref{tab:compatibilityClonage} features the results of the implementation of a normally distributed white-noise of standard-deviation of 0.1 and 0.0001 over two samples of 100 networks. Even with such small changes in the compatibility list, the number of edges drops significantly which dramatically alters the clustering coefficient making networks that cannot be considered as "similar" with the original one. One possible explanation can lie in the white-noise deteriorating the correlation between columns of the compatibility matrix. Moreover, this high sensitivity to small changes in the compatibility list should be seen with respect to the high robustness of the main benchmark characteristics of the original graph: when removing 5\% of the vertices the clustering coefficient follows a standard distribution $\mathcal{N}(0.5055,0.0148)$ (sample of 1000 trials, Shapiro test p-value of 0.04) whereas the original clustering coefficient was 0.5610. Similar results can be achieved with different methods: when removing 10\% of individuals and considering the resulting sub-graphs, we obtain an average global clustering coefficient of 0.5510 (standard deviation of 0.026 over 1000 observations), illustrating a high robustness of the network. Due to the high sensibility of the model to the compatibility list,  this method does not appear relevant. 

Our suggested methodology tries to take into account this criticism: the compatibility matrix has columns that are correlated because the compatibility list is formed by scalar-products of the individuals' characteristic vector. Hence, we suggest not directly altering the results of these scalar-products, but the the individuals' characteristic vector. 

Results are shown in Table \ref{tab:compatibilityClonage}. Rather than trying to infer a specific relationship between dispersion of the white-noise process and main benchmark characteristics (a relationship that would risk to be specific to this example), we only observe that the behavior of the output is much more continuous with the dispersion of the white-noise process than what we observed when we altered directly the compatibility list. Moreover, we find that, for this example, the dispersion $\sigma=0.1$ gives the best results, with networks that are very similar: the global clustering coefficient, the density, the average degree and the number of relations are very close to what characterizes the original graph. The networks obtained are also sufficiently different from the original to not to be considered as cloned: around 80\% of the edges that exist in the new networks did not exist in the original one.

\subsection{Decentralized optimization and agent-based model}

\begin{table*}[] 
\begin{adjustbox}{max width=\textwidth} 
  \begin{threeparttable}[t]
  \centering 
  \caption{\textbf{Random Effects regression estimates of the heterogeneity of the effect of the scope dimension ($\delta$) on the deviation of the equilibrium of the ABM} from the theoretical solution (SOP). In models (1) to (6), $\varkappa$ is the clustering coefficient of the optimal graph (solution of the SOP), in model (7) it is the Gini coefficient of the same graph \textit{(sample of 50 compatibility structures, simulated with $\iota=4$).}} 
  \label{tab:abmPanelDataScope2} 
\begin{tabular}{lccccccc} 
\midrule\midrule
 & \multicolumn{7}{c}{\textit{Dependent variable:}} \\ 
\cline{2-8} 
\\[-1.8ex] & Clustering & Distance & Avg. degree & Avg. Utility & S.d. Utility & Density & Gini\\ 
\\[-1.8ex] & (1) & (2) & (3) & (4) & (5) & (6) & (7)\\ 
\midrule
$\delta=0$ & $-$0.424$^{***}$ & $-$0.417$^{***}$ & $-$0.482$^{***}$ & $-$0.672$^{***}$ & 435,973 & $-$0.482$^{***}$ & $-$9.723$^{***}$ \\ 
  & {\scriptsize(0.067)} & {\scriptsize(0.076)} & {\scriptsize(0.023)} & {\scriptsize(0.014)} & {\scriptsize(392,293)} & {\scriptsize(0.023)} & {\scriptsize(1.571)} \\ 
 $\delta=1$ & 0.582$^{***}$ & 0.441$^{***}$ & $-$0.279$^{***}$ & $-$0.072$^{***}$ & 962,152$^{***}$ & $-$0.279$^{***}$ & 4.922$^{***}$\\ 
  & {\scriptsize(0.054)} & {\scriptsize(0.062)} & {\scriptsize(0.019)} & {\scriptsize(0.011)} & {\scriptsize(318,603)} & {\scriptsize(0.019)} & {\scriptsize(1.263)} \\ 
 $\delta=2$ & $-$0.207$^{***}$ & 0.084 & $-$0.301$^{***}$ & $-$0.103$^{***}$ & 791,552$^{**}$ & $-$0.301$^{***}$& 1.479 \\ 
  & {\scriptsize(0.054)} & {\scriptsize(0.062)} & {\scriptsize(0.019)} & {\scriptsize(0.011)} & {\scriptsize(318,603)} & {\scriptsize(0.019)} & {\scriptsize(1.263)} \\ 
 $\delta=3$ & $-$0.301$^{***}$ & 0.156$^{**}$ & $-$0.346$^{***}$ & $-$0.130$^{***}$ & 741,601$^{**}$ & $-$0.346$^{***}$ & 0.347\\ 
  & {\scriptsize(0.054)} & {\scriptsize(0.062)} & {\scriptsize(0.019)} & {\scriptsize(0.011)} & {\scriptsize(318,603)} & {\scriptsize(0.019)} & {\scriptsize(1.263)} \\ 
 $\delta=4$ & $-$0.322$^{***}$ & 0.249$^{***}$ & $-$0.364$^{***}$ & $-$0.141$^{***}$ & 899,695$^{**}$ & $-$0.364$^{***}$ & $-$3.700$^{**}$ \\ 
  & {\scriptsize(0.067)} & {\scriptsize(0.076)} & {\scriptsize(0.023)} & {\scriptsize(0.014)} & {\scriptsize(392,293)} & {\scriptsize(0.023)} & {\scriptsize(1.571)}\\ 
 \midrule
  $\varkappa$ & $-$0.011 & 0.468$^{*}$ & $-$0.152$^{**}$ & 0.154$^{***}$ & $-$647,334 & $-$0.152$^{**}$ & $-$53.950$^{***}$\\ 
  & {\scriptsize(0.113)} & {\scriptsize(0.278)} & {\scriptsize(0.077)} & {\scriptsize(0.017)} & {\scriptsize(1,692,293)} & {\scriptsize(0.077)} & {\scriptsize(9.179)} \\ 
 $\varkappa \mathbbm{1}_{[\delta=0]}$ & $-$0.318$^{**}$ & 0.765$^{***}$ & $-$0.254$^{***}$ & $-$0.325$^{***}$ & $-$335,625 & $-$0.254$^{***}$ & 53.950$^{***}$ \\ 
  & {\scriptsize(0.140)} & {\scriptsize(0.160)} & {\scriptsize(0.049)} & {\scriptsize(0.029)} & {\scriptsize(818,832)} & {\scriptsize(0.049)} & {\scriptsize(10.511)} \\ 
 $\varkappa \mathbbm{1}_{[\delta=1]}$ & $-$1.463$^{***}$ & $-$0.649$^{***}$ & $-$0.008 & $-$0.173$^{***}$ & $-$763,892 & $-$0.008 & $-$26.215$^{***}$\\ 
  & {\scriptsize(0.114)} & {\scriptsize(0.130)} & {\scriptsize(0.040)} & {\scriptsize(0.023)} & {\scriptsize(668,573)} & {\scriptsize(0.040)} & {\scriptsize(8.582)}\\ 
 $\varkappa \mathbbm{1}_{[\delta=2]}$& $-$0.196$^{*}$ & 0.144 & $-$0.060 & $-$0.192$^{***}$ & $-$640,581 & $-$0.060 & $-$8.289 \\ 
  & {\scriptsize(0.114)} & {\scriptsize(0.130)} & {\scriptsize(0.040)} & {\scriptsize(0.023)} & {\scriptsize(668,573)} & {\scriptsize(0.040)} & {\scriptsize(8.582)}\\ 
 $\varkappa \mathbbm{1}_{[\delta=3]}$ & $-$0.051 & 0.088 & 0.006 & $-$0.159$^{***}$ & $-$499,933 & 0.006 & $-$2.254 \\ 
  & {\scriptsize(0.114)} & {\scriptsize(0.130)} & {\scriptsize(0.040)} & {\scriptsize(0.023)} & {\scriptsize(668,573)} & {\scriptsize(0.040)} & {\scriptsize(8.582)} \\ 
 $\varkappa \mathbbm{1}_{[\delta=4]}$ & $-$0.023 & $-$0.043 & 0.046 & $-$0.146$^{***}$ & $-$644,439 & 0.046 & 22.598$^{**}$ \\ 
  & {\scriptsize(0.140)} & {\scriptsize(0.160)} & {\scriptsize(0.049)} & {\scriptsize(0.029)} & {\scriptsize(818,832)} & {\scriptsize(0.049)} & {\scriptsize(10.511)} \\ 
\midrule
 Cst. & 0.730$^{***}$ & 1.352$^{***}$ & 0.838$^{***}$ & 0.883$^{***}$ & 591,077 & 0.838$^{***}$ & 9.723$^{***}$ \\ 
  & {\scriptsize(0.054)} & {\scriptsize(0.133)} & {\scriptsize(0.037)} & {\scriptsize(0.008)} & {\scriptsize(807,459)} & {\scriptsize(0.037)} & {\scriptsize(1.378)} \\ 
\midrule
Adj. R$^{2}$ & 0.685 & 0.550 & 0.930 & 0.983 & 0.065 & 0.930 & 0.244\\ 
F Stat. & 1,098$^{***}$ & 623$^{***}$ & 6,680$^{***}$ & 29,537$^{***}$ & 47$^{***}$ & 6,680$^{***}$ &173$^{***}$ \\ 
\midrule \midrule
\textit{Note:}  & \multicolumn{7}{r}{500 observations, $^{*}$p$<$0.1; $^{**}$p$<$0.05; $^{***}$p$<$0.01} \\ 
\end{tabular}  

\end{threeparttable}
\end{adjustbox}
\end{table*}

We will now try to measure the distance between the decentralized equilibrium and the theoretical one (SOP). We first challenge the robustness of the generator to its input parameters, and then try to assess the effects of the frictions on the equilibrium.

First, despite the randomness of the generation seed, the intra-dispersion (the dispersion in results between the three runs associated with the same compatibility list) is negligible with respect to the dispersion observed between compatibility structures (always more than 10 times lower), with most indicators having dispersion of less than 2\% (max: 9\% in our sample) of their mean value. We observe a quadratic relation between the initialization density and the deviation from the theoretical graph for the clustering coefficient and the density, and long lasting effects on the output that need to be taken into account for later comparisons. Transitioning to a decentralized generation process reduces the maximum value reached by the objective function at equilibrium (optimality cost): the reduction of the average utility received by individuals is ranging from $6\%$ to $1\%$ for $\iota \in \{4,16,32,48\}$ (the seed density). This illustrates that the round-by-round decentralized process is able to find an equilibrium that is very close to the theoretical optimum.

Second, implementing scope limitations always reduces the utility reached by agents (objective function), while  decreasing the clustering coefficient, the average degree and the density (Table \ref{tab:abmPanelDataScope2}). The effect is significantly increasing with $\delta$, the scope depth. The utility dispersion is higher with scope limitations, but the impact of $\delta$ is unclear. Among all random effects models with clustered errors for the three runs, the average utility appears to be the most prominent and significant relationship, with very high significance of all regressors and a high adjusted $R^2$ (0.979), taking into account the effect of the initialization seed.

When controlling by the clustering coefficient (or the Gini coefficient) that has been observed in the SOP solution, clustering, distance, average utility and Gini coefficients are shown to be altered differently by the scope limitations, depending on the type of theoretical network. First, these estimates make the case $\delta = 0$ much more understandable: among all scope limitations, it is the most sensitive to the theoretical clustering of the graph, illustrating the inability of the agent-based model to generate highly clustered graphs with such scope limitations. Second, except for the Gini coefficient, the higher is $\delta$ (that is the lower are the scope limitations), the less are its effects dependent on the theoretical clustering coefficient, suggesting that this dependency would essentially matter for high scope limitations ($\delta \leq 1$). Last but not least, the optimality cost appears to be highly dependent of the theoretical optimal network clustering coefficient, with lower clustered networks being much more resilient to frictions. The cross-effect of scope limitations and theoretical clustering is also steadily decreasing with $\delta$. 

Overall, scope limitations have an ambiguous effect. As widely experimented in the Economics literature, imperfect information is reducing the quality of the optimum (the objective function reaches a lower value, here the total utility is lower). Nonetheless, this reduction of the quality of the equilibrium does not appear to be linearly decreasing with the amount of information accessible to agents, suggesting a quadratic relation with information. The network benchmark properties are affected heterogeneously, but all resulting networks share a common reduction of the total number of edges and nodes' degree, inducing a reduction of the global density, an increase in the mean distance and an alteration of the global clustering coefficient. All compatibility structures exhibits similar reactions to frictions, but less clustered networks have proved slightly more resilient, especially when frictions are very high.

\section{Discussion} \label{sec:discussion}
In the sub-field of graph reproduction, our results with the Optimal Graph generator (from the SOP, Definition \ref{def_SOP}) have illustrated that not only it is possible to use our objective-maximization framework to model human networks, but also to generate similar networks from an original one \cite{ali2014synthetic, nettleton2016synthetic, johnson2024epidemic}. This suggests that it is perfectly possible to conciliate reproducibility and fidelity with the original graph, with interpretability of the generation process. 

Moreover, our micro-based framework could be extended to numerically estimate the utility function parameters or compatibilities between individuals using real graphs. This would help link pure-graph features (\textit{e.g.}, clustering coefficient, degree distribution) to the behavior of network agents. However, this requires addressing two degrees of freedom: the compatibility matrix and the utility function parameters. Two approaches could be explored: inferring parameters from a compatibility distribution or inferring compatibilities from fixed parameters. The first approach standardizes compatibilities by settings the parameters of the utility function, allowing comparisons. However, as determining these distributions and effectively setting the utility functions seems complex, the second approach may be preferred. It involves fixing parameters to derive compatibilities, enabling comparisons if graphs are supposed to share the same utility parameters.

The question of utility function selection is important for several reasons. First even if we can generalize our results for a class of functions (cf. supplementary materials), one can argue that we have no guarantee that the agent-based model would present similar properties with a different function. If our function is more complex than what was previously experimented \cite{jackson2007study, hellmann2014} making relationships both substitutable and not linear in their costs, several other determinants could be attached. The centrality implemented in \cite{jackson1996} could be added by making the advantage of having a relationship with two individuals depend on the fact that these two individuals are themselves related. In addition, the choice of the utility function raises the question of the computational efficiency of both generators. The downside of our fully micro-based approach is that it is comparatively less efficient than other generators that do not care for consistency of agents behaviors at each scale. Some later developments could therefore explore new utility functions that could induce computational efficiency gains, allowing the generator to handle higher dimensional networks. 

Furthermore, this work strengthens the need for further developments on the consequences of the latent space geometry. Several works including \cite{saha_study_2024, papamichalis_latent_2022} have tried to show that the degree distribution can be a power law distribution only using non-euclidean geometry, but the theoretical implications of using a type of geometry rather than another remain an open question today. Our choice of using a scalar product as a compatibility measure can be seen as implementing a non-euclidean latent space with $d^{\text{latent}}(x,y) = \left(\chi(x,y)\right)^{-1}$ (for $x\neq y$, 0 elsewhere). This conic non-linear geometry has allowed us to recover similar results as the one found with hyperbolic geometry, without clear explanation. 

Finally, our approach paves the way for studying the consequences of real world imperfections on the shape of networks with our generators. Indeed, the Optimal Graph generator can be used as a proxy for what would have been obtained with independent and rational agents in conditions given by Theorem \ref{thm_invisibleHand}. On the other hand, the agent-based model can be used to obtain estimates of the impact of real-life constrains such as geographic or normative ones.

\section*{Conclusion} \label{sec:conclusion}
Overall, our main contribution has been to introduce a mathematical framework able to conciliate three distinct objectives: fidelity to empirical network when generating from real-world data, consistency at every scale in the generation process with fully individual-based postulates, and versatility in the implementation allowing for the use of agent-based round-by-round simulations. Indeed, our Optimal graph generator has been proven to be able to regenerate any real-life network (that satisfy very loose criteria) without leaving aside our agent-based assumptions. Moreover, we have illustrated on an example that our framework was even able to generate very similar artificial networks from an empirical one.

Our micro-based consistency has enabled us to build an agent-based simulator with round-by-round rules, with approximately zero optimality cost, justifying the relevance of the coexistence of the two models: if the agent-based model allows to test hypothesis and real-world constrains, the analytical Optimal graph model offers a quicker and a more versatile generator for users that would not care for such hypothesis tests.  This has enabled us to challenge the effect of two real-world constrains, time and limited information, on the resulted network: for the limited information constrain, highly clustered networks and populations with highly dispersed affinity scores have proven to be less resilient. 

\section*{Code and data availability}
Our implementation of the generators is available in C++ on the \href{https://github.com/Aldric-L/SocialMatrixSimulation}{GitHub repository} where both the source code and the executable can be found under the \href{https://www.gnu.org/licenses/gpl-3.0.html}{GNU General Public License v3.0}. All simulations data are also provided, with the R script that enabled us to create all tables and to compute all statistical estimations.

%\nocite{*}

\bibliography{main}

\clearpage
\onecolumngrid
\appendix
\section{Theoretical and mathematical Framework}

\subsection{Fundamental definitions}

\begin{defi}[Utility function]\label{app_def_utilityFunc}
For a given $I \in \mathbb{N}^*$, a utility function $U(\alpha, C)$ is an element of $\mathcal{U}_I$, the set of applications $U : \mathcal{R}^{I-1} \times (\mathbb{R}^{*}_{+})^{I-1} \longrightarrow \mathbb{R}$ with $[0 ,1) \subseteq \mathcal{R} \subseteq \mathbb{R}$, and $\mathcal{R}^{I-1}$ convex. Every element $U : (\alpha, C) \longmapsto \mathbb{R}$ of $\mathcal{U}_I$ verifies:
\begin{enumerate}
\item $\forall C \in (\mathbb{R}^{*}_{+})^{I-1}, U(0 , C) = 0 $, and $\frac{\partial U}{\partial \alpha_i}(0 ) > 0 $
\item with $H$ bijective and increasing in its second component,\begin{equation*} \exists \begin{cases}H : \mathbb{R}^{*}_{+} \times \mathcal{R} \longrightarrow \mathbb{R}, \quad H \in \mathscr{C}^1\\ S : \mathcal{R}^{I-1} \longrightarrow \mathbb{R}, \quad S \in \mathscr{C}^1\end{cases} \quad \text{s.t. }\forall i, \frac{\partial U}{\partial \alpha_i} = H(c_i, \alpha_i) + S(\alpha)
\end{equation*}

\item \begin{equation*} \exists \begin{cases}h : \mathbb{R}^{*}_{+} \times \mathcal{R} \longrightarrow \mathbb{R}, \quad h \in \mathscr{C}^0 \\ s : \mathcal{R}^{I-1} \longrightarrow \mathbb{R}, \quad s \in \mathscr{C}^0 \end{cases} \text{ s.t. } \forall j \neq k, \quad \frac{\partial^2 U}{\partial \alpha_j \alpha_k} = s(\alpha), \quad \frac{\partial^2 U}{\partial \alpha_j^2} = h(c_j, \alpha_j) + s(\alpha)\end{equation*}
with $\begin{cases}S(\alpha) = 0  \vee s(\alpha) = 0  \Leftrightarrow \alpha = 0  \\ S(\alpha) = 0  \vee s(\alpha) = 0  \Rightarrow \forall j \leq I-1, h(c_j, \alpha_j) = 0  \wedge H(c_j, \alpha_j) = 0 \\
\forall \mathcal{A} \in \mathcal{R}^{I-1}, \exists \bar s \in \mathbb{R}, \forall a \in \mathscr{S}(\mathcal{A}), S(a) = \bar s \\
\forall \alpha, \forall i, h(\alpha_i, c_i) \leq 0 , s(\alpha) \leq 0 , S(\alpha) \leq 0 \end{cases}$ 
\end{enumerate}
\end{defi}

\begin{lem}\label{app_lem_utilityAreConcave}
All Utility functions are strictly concave in their first argument on $\mathcal{R}^{I-1}$.
\end{lem}

\begin{proof}[Proof (Lemma \ref{app_lem_utilityAreConcave})]
Let's denote the Hessian of $U \in \mathcal{U}$ and take $(C, \alpha) \in (\mathbb{R}^{*}_{+})^{I-1} \times \mathcal{R}^{I-1}$ with $\alpha \neq 0 $,
\begin{align*}
&\mathbf{H} =
\begin{pmatrix} 
h(c_1, \alpha_1) + s(\alpha) & s(\alpha) &\hdots &  s(\alpha) \\
s(\alpha) & h(c_2, \alpha_2) +  s(\alpha)& \ddots &\vdots  \\
\vdots & \ddots & \ddots &s(\alpha)  \\
s(\alpha) & \hdots & s(\alpha) & h(c_{I-1}, \alpha_{I-1}) + s(\alpha)  \\
\end{pmatrix} \\
&\Rightarrow \text{det}\left(\mathbf{H}\right) = \begin{vmatrix} 
h(c_1, \alpha_1) + s(\alpha) & s(\alpha) &\hdots &  s(\alpha) \\
s(\alpha) & h(c_2, \alpha_2) +  s(\alpha)& \ddots &\vdots  \\
\vdots & \ddots & \ddots &s(\alpha)  \\
s(\alpha) & \hdots & s(\alpha) & h(c_{I-1}, \alpha_{I-1}) + s(\alpha)  \\
\end{vmatrix}
\\
&\Rightarrow \text{det}\left(\mathbf{H}\right) = \begin{vmatrix} 
h(c_1, \alpha_1)  & 0  &\hdots & 0  \\
0  & h(c_2, \alpha_2) & \ddots &\vdots  \\
\vdots & \ddots & \ddots &0   \\
0  & \hdots & 0  & h(c_{I-1}, \alpha_{I-1})  \\
\end{vmatrix}
\end{align*}

As $h$ and $s$ have been assumed to be strictly negative for $\alpha \neq 0 $, we have that $\mathbf{H}$ is negative definite on $\mathcal{R}^{I-1} \backslash \{0 \}$. Let's prove strict concavity over $\mathcal{R}^{I-1}$. We should have $\forall (x,y) \in (\mathcal{R}^{I-1})^2, x \neq y, \forall \lambda \in (0 ,1), U(\lambda x + (1-\lambda)y, C) < \lambda U(x, C) + (1-\lambda) U(y, C)$.
We can distinguish two cases:
\begin{itemize}
\item $x\neq0, y\neq 0 $, the standard second characterization of strict concavity provides the proof as $\mathbf{H}$ is negative definite on the domain of $x$ and $y$.
\item $x=0 \text{ or }y=0 $ (let's take $x=0 $): $U(\lambda x + (1-\lambda)y, C) = U((1-\lambda)y, C) $. As we have strict concavity everywhere but in 0  and as $U(0 , C) = 0 $ we have 
$$U((1-\lambda)y, C) < (1-\lambda)U(y, C) \Rightarrow U(\lambda x + (1-\lambda)y, C) < (1-\lambda)U(y, C) =  \lambda U(x, C) + (1-\lambda)U(y, C)$$
\end{itemize}
\end{proof}

\begin{defi}[Individuals space]\label{app_def_indivSpace}
We denote as \textbf{Individuals space} the product space $$\mathcal{I} = (\mathcal{E}, \innerproduct{.}{.}) \times \mathcal{U}_I$$ with $(\mathcal{E}, \innerproduct{.}{.})$ a finite-dimensional euclidean space defined on the vector space $\mathcal{E}$ with $\innerproduct{.}{.}$ an inner product over $\mathbb{R}^{+}$. Individuals sets are elements of $\mathcal{I}^I$. No individual should be orthogonal to any other in its first component.
\end{defi}

\begin{defi}[Individuals graph]\label{app_def_indivGraph}
For a given individuals space $\mathcal{I}$ and a given individuals set $\mathscr{I} \in \mathcal{I}^I$, an \textbf{individual graph} is a 3-tuple $\mathcal{G}(\mathscr{I} , \alpha, C)$ with $\alpha, C$ two vectors of $(\mathbb{R}^{*}_{+})^{n_I}$ and $(\mathbb{R}_{+})^{n_I}$) with $n_I=\binom{I}{2} = \frac{I^2-I}{2}$ :
\begin{itemize}
    \item $\alpha$ is a \textbf{weight list}. $$\alpha = (\alpha_{1,2}, \hdots, \alpha_{1,I},  \alpha_{2,3}, \hdots, \alpha_{2,I}, \hdots,\alpha_{I-1,I} )^T$$
    with $\alpha_{i,j}$ the edge weight between agent i and j
    \item $\chi(i,j) = \innerproduct{P_i}{P_j}\mathbbm{1}(i\neq j)$ is the \textbf{compatibility function}, with $(i,j)\in \mathscr{I}^2$ and $P_i$ the first component of i (in $\mathcal{E}$).
    \item $C$ is a \textbf{compatibility list} $$C = (c_{1,2}, \hdots, c_{1,I},  c_{2,3}, \hdots, c_{2,I}, \hdots,c_{I-1,I} )^T$$ with $c_{i,j} = \chi(i,j) = \innerproduct{P_i}{P_j}$.
    
\end{itemize}

On each graph is defined $d : \mathcal{I}^2 \rightarrow \mathbb{N}$ a distance between agents $i$ and $j$: the smallest number of intermediary non-null edges it takes to reach a vertex from another.
\end{defi}

\begin{rmk}\label{app_rmk_indivUtility}
In a Individuals graph, the individual utility of agent $i$, $U_I \in \mathcal{U}_I$ is evaluated on $\alpha^i$ and on $C^i$, that is on the vectors of weights and compatibilities with respect to agent i.
\end{rmk}

\begin{nota}\label{app_nota_cAndAlpha}
We denote $\forall i \in \Iintv*{1, I }, \alpha^i = (((\alpha_{i,j})_{i<j\leq I}), ((\alpha_{j,i})_{1\leq j<i}))$ the vector of all vertices liked with the vertex $i$ of dimension $I-1$. %$\alpha_i$ is the ith-component of $\alpha$, and 
$\alpha^i_j$ is the jth-component of $\alpha^i$. When we write $\alpha_{i,j}$, we consider $\alpha_{\min(i,j), \max(i,j)}$. We use the same notations for $C$ (with $C^i$ and $c^i_j$).
\end{nota}

\subsection{A network optimization problem}

\begin{defi}[Individual Optimization Problem]\label{app_def_IOP}
For a given individuals graph $\mathcal{G}_N(\mathscr{I} , \alpha, C)$, we define for all individual $i \in \Iintv*{1, N }$ the Individual Optimization Problem (IOP) as:
\begin{equation}
\begin{cases}
\argmax_{\alpha^i} U_i(\alpha^i, C^i) \\
\text{s.t.} \quad \alpha^i \in [0 ,1)^{I-1}
\end{cases}
\end{equation}
\end{defi}

\begin{thm}[Unicity of the individual optimal graph]\label{app_thm_unicityIOP}
For any $\mathcal{G}(\mathscr{I} , \alpha, C)$, and for all individual $i \in \Iintv*{1, N }$ the Individual Optimization Problem (IOP) has only one unique solution.
\end{thm}

\begin{proof}[Proof (Theorem \ref{app_thm_unicityIOP})]
By Lemma \ref{app_lem_utilityAreConcave}, we have that $U_i$ is strictly concave in its first argument on $\mathcal{R}^{I-1}$, and on $[0 ,1)^{I-1}$ as it is a convex subset of $\mathcal{R}^{I-1}$. As constraints are qualified and convex, we have that the IOP has one unique solution. 
\end{proof}

\begin{defi}[Social Optimization Problem]\label{app_def_SOP}
For a given individuals graph $\mathcal{G}_N(\mathscr{I} , \alpha, C)$, we define the Social Optimization Problem (SOP) as:
\begin{equation}
\begin{cases}
\argmax_{\alpha} \mathcal{O} (\alpha, C) =  \sum^N_{i=1} U_i(\alpha^i, C^i)\\
\text{s.t} \quad \begin{cases} 
\alpha \in [0 ,1)^{n_N} \\
\nexists i \in \Iintv*{1, N }, \forall j \in \Iintv*{1, N-1 }, \alpha^i_j = 0 
\end{cases}
\end{cases}
\end{equation}
{\footnotesize With $U_i$ the utility function of agent $i$, that is its second component in the individuals space.}
\end{defi}

\begin{lem}[Strict concavity of the objective]\label{app_lem_concavityObj}
For a given individuals graph $\mathcal{G}_I(\mathscr{I} , \alpha, C)$,  $\mathcal{O}(C, \alpha) =  \sum^I_{i=1} U_i(C^i, \alpha^i)$ is strictly concave on the feasible set $\{\alpha \in \mathcal{R}^{I-1}, [\alpha \in [0 ,1)^{n_I}] \wedge [\nexists i \in \Iintv*{1, I }, \forall j \in \Iintv*{1, I-1 }, \alpha^i_j = 0 ] \}$.
\end{lem}

\begin{proof}[Proof (Lemma \ref{app_lem_concavityObj})]
Let's take $\alpha$ in $\mathcal{R}^{I-1}$ and $\forall (i,j,m,p) \in (\Iintv*{1, I })^4, i\neq j\neq m \neq p$:
\begin{align*}
\frac{\partial \mathcal O}{\partial \alpha_{i, j}} = \frac{\partial U_i}{\partial \alpha_{i, j}} + \frac{\partial U_j}{\partial \alpha_{i, j}}
\Rightarrow
\begin{cases}
\frac{\partial^2 \mathcal O}{\partial \alpha_{i, j}^2} = h_i(C^i, \alpha^i) + s_i(\alpha^i) + h_j(C^j, \alpha^j) + s_j(\alpha^j) \\
\frac{\partial^2 \mathcal O}{\partial \alpha_{i, m} \partial \alpha_{i, j}} = s_i(\alpha^i)  \\
\frac{\partial^2 \mathcal O}{\partial \alpha_{p, m} \partial \alpha_{i, j}} = 0 
\end{cases}
\end{align*}

Let's now prove that the Hessian matrix ($\mathbf{H}$) of $\mathcal O$ is negative definite. We will, from now on consider only $-\mathbf{H}$ and prove that it is definite positive, under the assumption that $\nexists (i,j) \in (\Iintv*{1, I })^2 / i\neq j, s_i(\alpha^i) = 0  \wedge s_j(\alpha^j) = 0  $. We prove it by induction over $I$ ($I>1$). We denote $\mathbf{H}_I$ the Hessian matrix for I, $\mathbf{H}_{I,(i,j)}$ its ijth-component, and $n_I = \binom{I}{2}$ the dimension of $\alpha$ and $C$.

\underline{Initialization:} $I=2 \Rightarrow n_2=1$ $\text{det}(-\mathbf{H}_2) = - (h_1 + s_1 + h_2 + s_2) > 0 $.
 
\underline{Inheritance:} Assume for a given $I>1$ that $-\mathbf{H}_I$ is positive definite. 
$$\mathbf{H}_{I+1} = \begin{pmatrix} \mathbf{H}_{I} & R \\ R^T & Q\end{pmatrix}$$
with $\mathbf{H}_{I+1} \in \mathcal{M}_{n_{I+1}}(\mathbb{R}), \mathbf{H}_{I} \in \mathcal{M}_{n_{I}}(\mathbb{R}), Q \in \mathcal{M}_{(n_{I+1}-n_{I})}(\mathbb{R}), R \in \mathcal{M}_{n_{I}\times (n_{I+1}-n_{I})}(\mathbb{R})$, $\mathbf{H}$ and Q being both symmetric matrices by construction.
$$R = \begin{pmatrix} 
\frac{\partial^2 \mathcal O}{\partial \alpha_{1, 2} \partial \alpha_{1, I+1}} &\hdots &  \frac{\partial^2 \mathcal O}{\partial \alpha_{1, 2} \partial \alpha_{I, I+1}} \\
\vdots & &\vdots \\
\frac{\partial^2 \mathcal O}{\partial \alpha_{I-1, I} \partial \alpha_{1, I+1}} & \hdots & \frac{\partial^2 \mathcal O}{\partial \alpha_{I-1, I} \partial \alpha_{I, I+1}} 
\end{pmatrix}
\quad 
Q = \begin{pmatrix} 
\frac{\partial^2 \mathcal O}{\partial \alpha_{1, I+1}^2} &\hdots &  \frac{\partial^2 \mathcal O}{\partial \alpha_{1, I+1} \partial \alpha_{I, I+1}} \\
\vdots & \ddots &\vdots \\
\frac{\partial^2 \mathcal O}{\partial \alpha_{I, I+1} \partial \alpha_{1, I+1}} & \hdots & \frac{\partial^2 \mathcal O}{\partial \alpha_{I, I+1}^2}  
\end{pmatrix}
$$

To prove that $-\mathbf{H}_{I+1}$ is positive definite, we proceed by leading principal minors. As $-\mathbf{H}_{I}$ is assumed to be positive definite, we only need to consider the $n_{I+1}-n_{I}$ leading principal minors. We take $0  < i \leq n_{I+1}-n_{I}$ and denote $q_{i,j}$ the elements of $Q$, $R_i$ the ith column of $R$ and $\mathbf{H}_{I+1, n_I+i}$ the matrix formed by the first $n_I+i$ rows and columns of $\mathbf{H}_{I+1}$. The ith leading minor ($\mu_i$) that we have to consider can be written as:
$$\mu_i = \vert -\mathbf{H}_{I+1, n_I+i} \vert = \begin{vmatrix} -\mathbf{H}_{I+1, n_I+i-1} & -(R_i^T,q_{1,i},\hdots, q_{i-1,i})^T  \\ -(R_i^T,q_{1,i},\hdots, q_{i-1,i}) & -q_{i,i}\end{vmatrix}$$. 

As $R_i$ is the vector $(\frac{\partial^2 \mathcal O}{\partial \alpha_{1, 2} \partial \alpha_{i, I+1}}, \hdots, \frac{\partial^2 \mathcal O}{\partial \alpha_{I-1, I} \partial \alpha_{i, I+1}})^T$, which is composed of at least one $s_i$ (the others elements being 0 ), and $(q_{1,i},\hdots, q_{i-1,i})$ being $(s_{I+1},\hdots,s_{I+1})$, we have that $(R_i^T,q_{1,i},\hdots, q_{i-1,i}) \neq 0_{\mathbb{R}^{n_I+i}}$ iff $s_i$ (or $s_{I+1}$ if $i>1$) is non-null. 

\begin{itemize}
\item if $s_i$ (and $s_{I+1}$ if $i>1$) is null:
$$\mu_i = \text{det}(-q_{i,i})\text{det}\left( (-\mathbf{H}_{I+1, n_I+i-1}) \right) > 0 $$

\item If $s_i$ (or $s_{I+1}$ if $i>1$) is non-null, we have that $-(R_i^T,q_{1,i},\hdots, q_{i-1,i})$ is positive and $-q_{i,i}$ is strictly non-negative, and as $-\mathbf{H}_{I}$ is definite positive (and so is $(-\mathbf{H}_{I})^{-1}$), we have that $\mu_i > 0 $ by the block matrix determinant formula: 
$$\mu_i = \text{det}(-q_{i,i})\text{det}\left( (-\mathbf{H}_{I+1, n_I+i-1}) -\frac{1}{q_{i,i}} (-(R_i^T,q_{1,i},\hdots, q_{i-1,i})( -(R_i^T,q_{1,i},\hdots, q_{i-1,i})^T )) \right)$$
\begin{align*}
\mu_i > 0  \Leftrightarrow& \text{ det}\left( (-\mathbf{H}_{I+1, n_I+i-1}) -\frac{1}{q_{i,i}} (R_i^T,q_{1,i},\hdots, q_{i-1,i})(R_i^T,q_{1,i},\hdots, q_{i-1,i})^T \right) >0 \\
\Leftrightarrow& \text{ det}\left( (-\mathbf{H}_{I+1, n_I+i-1}) -\frac{1}{q_{i,i}} (R_i^T,q_{1,i},\hdots, q_{i-1,i})(R_i^T,q_{1,i},\hdots, q_{i-1,i})^T  \right)\\
& \geq \text{det}\left(-\mathbf{H}_{I+1, n_I+i-1}\right) \underbrace{-\frac{1}{q_{i,i}} \text{det}\left((R_i^T,q_{1,i},\hdots, q_{i-1,i})(R_i^T,q_{1,i},\hdots, q_{i-1,i})^T \right)}_{=0 } > 0  
\end{align*}
\end{itemize}

Since $-\mathbf{H}_{I+1, n_I+i-1}$ was definite positive (if $i=1, \mathbf{H}_{I+1, n_I+i-1}$ is $\mathbf{H}_{I}$), and as $\mu_i>0 $, $-\mathbf{H}_{I+1, n_I+i}$ is definite positive. We repeat this procedure until $i=n_{I+1}-n_{I}$ and get that $-\mathbf{H}_{I+1}$ is definite positive too.

If the assumption that $\nexists (i,j) \in (\Iintv*{1, I })^2, i\neq j, s_i(\alpha^i) = 0  \wedge s_j(\alpha^j) = 0  $ does not hold, it is easy to notice that $\mathbf{H}_I$ will have at least one column and row that is full of 0 , thus $\text{det}(\mathbf{H}_I) = 0 $ in this case.

Hence, in the feasible set, $\mathbf{H}_{I}$ is negative definite and as the feasible set is trivially connected we have that $\mathcal{O}$ is strictly concave on the feasible set.

\end{proof}

\begin{thm}[Unicity of the optimal graph]\label{app_thm_unicity}
For any $\mathcal{G}(\mathscr{I} , \alpha, C)$, the Social Optimization Problem has only one unique solution.
\end{thm}

\begin{proof}[Proof (Theorem \ref{app_thm_unicity})]
We can define the feasible set as $\{\alpha \in \mathcal{R}^{I-1}, [\alpha \in [0 ,1)^{n_I}] \wedge [\nexists i \in \Iintv*{1, I }, \forall j \in \Iintv*{1, I-1 }, \alpha^i_j = 0 ] \}$, and observe that it is convex. By Lemma \ref{app_lem_concavityObj}, we have that the objective function (restricted to its first argument) is strictly concave on the feasible set. As constraints are qualified and convex, we have that the SOP has one unique solution. 
\end{proof}

\begin{rmk}
We denote as "unconstrained solution" for a constrained problem a solution that corresponds to a critic point of the objective function.
\end{rmk}

\begin{thm}[Marginal utilities]\label{app_thm_marginalUtilities}
For a given individuals graph $\mathcal{G}(\mathscr{I} , \alpha, C)$, if $\alpha$ the solution of the Social Optimization Problem is unconstrained, the marginal utility of agent $i$ with respect to $\alpha_{i, j}$ is equal to the reciprocal of the marginal utility of agent $j$ with respect to $\alpha_{i,j}$.
\end{thm}

\begin{proof}[Proof (Theorem \ref{app_thm_marginalUtilities})]
For $i$ and $j$ different and in $\Iintv*{1, I }$, $\forall k \in \{i,j\}, \exists m_k \in \Iintv*{1, I-1 }, \alpha^k_{p_k}=\alpha_{i,j}$ by construction of $\alpha$.
\begin{equation*}
\frac{\partial \mathcal O}{\partial \alpha_{i, j}} = \sum^I_{i=1} \frac{\partial U_i}{\partial \alpha_{i, j}} = 0  \Leftrightarrow  \frac{\partial U_i}{\partial \alpha^i_{p_i}} + \frac{\partial U_j}{\partial \alpha^j_{p_j}} = 0  \Leftrightarrow \frac{\partial U_i}{\partial \alpha_{i, j}} = -\frac{\partial U_j}{\partial \alpha_{i, j}}
\end{equation*}
\end{proof}

\begin{thm}[Invisible hand]\label{app_thm_invisibleHand}
For any $\mathcal{G}(\mathscr{I} , \alpha, C)$ with $\mathcal{C}$ the compatibility matrix defined as $\mathcal{C} = (c_{i,j})_{1\leq i,j\leq N}$ ($\mathcal{C}^{i}$ its ith-column), and $(\alpha^i)^*$ the solution vector of the Individual optimization problem (IOP) evaluated on $c^i$, the compatibility list derived from $\mathcal{C}_{i}$, if
\begin{enumerate}
\item $\exists \bar c \in \mathbb{R}^{N} \text{s.t. } \forall i \in \Iintv*{1, N }, \mathcal{C}^i \in \mathscr{S}(\bar c)$, with $\mathscr{S}(\bar c)$ the symmetric group of $\bar c$
\item there exists $(H,S)$ two functions satisfying conditions of Definition \ref{app_def_utilityFunc} such that $\forall i \in \Iintv*{1, N }, (\mathscr{I})_i = (P_i, U_i), \forall j \in \Iintv*{1, N-1 }, \frac{\partial U_i}{\partial \alpha^i_j} = H(c_j^i, \alpha_j^i) + S(\alpha^i)$
\end{enumerate}
Then there exists $\hat \alpha  \in \mathbb{R}^{n_N}$ such that $\forall i \in \Iintv*{1, N }, (\alpha^i)^* \in \mathscr{S}(\hat \alpha)$. Moreover, $\alpha^{**}$ the solution of the associated SOP, verifies: $$\forall (i,j) \in \Iintv*{1, N }^2, \alpha^{**}_{i,j} = (\alpha^{\max{\{i,j\}}})^*_{\min{\{i,j\}}}$$
 \end{thm}

\begin{proof}[Proof (Theorem \ref{app_thm_invisibleHand})] 
Condition (2) implies that there exists some utility function $\bar U$ such that $\forall i \in \Iintv*{1, N }, U_i = \bar U + \epsilon_i$ with $\epsilon_i$ some real constant, \textit{i.e.} all Individuals Optimization Problem (IOP) are the same. Condition (1) implies that each solution $(\alpha^i)^*$ (which exists by Theorem \ref{app_thm_unicityIOP}) must be a permutation of a vector $\hat \alpha$ of the same domain of $\alpha^i$.  As compatibilities are symmetric, even if the $(\alpha^i)^*$ are permutations, we have that the weight of the vertex that links individual $i$ and $j$ will be the same in both $(\alpha^i)^*$ and $(\alpha^j)^*$.

We prove this latter statement by contradiction. Let's take $C$ that satisfies condition (1) and consider ($(\alpha^k)^*, (\alpha^m)^*$) the solution vectors of the IOP associated with two compatibility lists. For the sake of contradiction, we assume that $ \alpha^k_m \neq \alpha^m_k$ but $((\alpha^k)^*, (\alpha^m)^*) \in (\mathscr{S}(\hat \alpha))^2$. We would have $ \frac{\partial U_m}{\partial \alpha^m_k}(a) = H(c_k^m, a_k) + S(a),\frac{\partial U_k}{\partial \alpha^k_m}(a)= H(c_m^k, a_m) + S(a)$. By assumption, we have that $\alpha^k_m \neq \alpha^m_k$, implying that $\exists a \in \mathscr{S}(\hat \alpha), \frac{\partial U_m}{\partial \alpha^m_k}(a)\neq\frac{\partial U_k}{\partial \alpha^k_m}(a)$ but as $c_k^m=c_m^k$ and $\forall a \in \mathscr{S}(\hat \alpha), S(a) = \bar s \in \mathbb{R}$, the two functions are the same by (2), which is absurd.

Since we know by Theorem \ref{app_thm_unicity} that there exists a decentralized solution, and as the objective function of the SOP is $\mathcal{O} = \sum_{i\in \Iintv*{1, N }}U_i$, if $\alpha^{**}$ maximizes each $U_i$, it is the maximizer of $\mathcal{O}$. By definition $U_i$ only depends on: $$(\alpha^{**})^i = ((\alpha^{\max{\{i,j\}}})^*_{\min{\{i,j\}}})_{1\leq j\leq N, i\neq j}^i = ((\alpha^{i})^*_{1}, ..., (\alpha^{i})^*_{i-1}, (\alpha^{i+1})^*_{i}, ..., (\alpha^{N})^*_{i} )$$
By the symmetry property that we have previously established, it can be rewritten as $(\alpha^{**})^i =((\alpha^{i})^*_{1}, ..., (\alpha^{i})^*_{i-1}, (\alpha^{i})^*_{i+1}, ..., (\alpha^{i})^*_{N} )$ which is by assumption a maximizer of $U_i$.

\end{proof}

\begin{rmk}
Theorem \ref{app_thm_invisibleHand} states that, if the two conditions are verified, all pairs of individuals $(i,j)$ will desire the same $\alpha_{i,j}$ in the two individuals optimization problems of agents i and j (IOPs), and in the whole SOP.
\end{rmk}

\begin{thm}\label{app_thm_alphaIsAlwaysOpti}
For any individuals set $\mathscr{I} \in (\mathcal{E}, \innerproduct{.}{.}) \times \mathcal{U}_N$ , and weight list $\alpha$ of dimension $n_N = \binom{N}{2}$ satisfying:
\begin{itemize}
\item $\forall a \in \alpha, a \in [0 ,1)$
\item $\nexists i \in \Iintv*{1, N }, \forall j \in \Iintv*{1, N-1 }, \alpha^i_j = 0 $
\item $\forall M \leq I, \left(\sum_{i\leq M} H_i\right)\left( C, \alpha \right)$ is bijective in its first component.
\end{itemize}

There exists:
\begin{enumerate}
\item a unique compatibility list $\bar C$ such that $\alpha$ is the unique unconstrained solution of the Social Optimization Problem for the individuals graph $\mathcal{G}(\mathscr{I} , \alpha, \bar C)$
\item a set of compatibility lists $\underbar C$ such that $\forall C \in \underbar C$, $\alpha$ is the unique solution of the Social Optimization Problem for the individuals graph $\mathcal{G}(\mathscr{I} , \alpha, C)$. 
$$\underbar C = \left\{ C \in (\mathbb{R}^{*}_{+})^{n_I}, \quad c_{ij} \in \begin{cases} \{\bar c_{ij}\} & \text{if } \alpha_{ij} > 0 \\  (0 , \bar c_{ij}]& \text{if } \alpha_{ij} = 0 \end{cases}\right\}$$
\end{enumerate}
\end{thm}

\begin{proof}[Proof (Theorem \ref{app_thm_alphaIsAlwaysOpti})]
Let's take $\mathscr{I} \in (\mathcal{E}, \innerproduct{.}{.}) \times \mathcal{U}_N$.

For point (1):  We can construct a unique compatibility list ($C$) without any condition on $\alpha$ (for $i<j\leq N$).
\begin{align*}
(\nabla \mathcal{O})(\alpha) = 0  \Leftrightarrow \begin{cases} \vdots \\ \frac{\partial U_i}{\partial \alpha_{i,j}} + \frac{\partial U_j}{\partial \alpha_{i,j}} = 0  \\ \vdots \end{cases}
\Leftrightarrow  \begin{cases} \vdots \\ H_i(c_{i,j}, \alpha_{i,j}) + H_j(c_{i,j}, \alpha_{i,j}) + s_i(\alpha^i) + s_j(\alpha^j) = 0  \\ \vdots \end{cases}
\end{align*}

As $\alpha$ is taken here as given (we are solving for $C$), each $c_{i,j}$ is determined by a single equation. As we have assumed that $\forall (i,j) \in (\Iintv*{1, N})^2, H_i+H_j$ is bijective, we have that there exists a single $C$ verifying $(\nabla \mathcal{O})(\alpha) = 0 $ for $\alpha$ (and the individuals set) given.

From this compatibility list ($C$) and this set of individuals, we can form an individuals graph and write its associated SOP (Definition \ref{app_def_SOP}). By Theorem \ref{app_thm_unicity} this SOP has only one solution. By construction of $C$, $\alpha$ satisfies the first order condition of the optimization problem, and by assumption, $\alpha$ is in the feasible set. Hence, $\alpha$ is the unique solution of an SOP for the $\mathcal{G}(\mathscr{I} , \alpha, C)$ we have constructed.

For point (2): We can write the lagrangean of the Social Optimization problem (Definition \ref{app_def_SOP}):
\begin{equation*}
\mathcal L (\alpha, C, \lambda, \mu) = -  \mathcal{O}(\alpha, C) - \sum_{k=1}^{n_I} \lambda_k \alpha_k - \sum_{k=1}^{n_I} \mu_k (1-\alpha_k)
\end{equation*}

We have already proven that the optimal solution of this problem is unique. We can write the first order conditions that characterize this solution:
\begin{equation*}
\frac{\partial \mathcal L}{\partial \alpha_{i,j}} = - H_i(c_{i,j}, \alpha_{i,j}) - H_j(c_{i,j}, \alpha_{i,j}) - S(\alpha^i) - S(\alpha^j) - \lambda_{i,j} + \mu_{i,j} = 0 
\end{equation*}

As $\forall (i,j) \in (\Iintv*{1, N})^2, H_i+H_j$ is bijective, solving for $c_{i,j}$ is possible and has only one solution:
\begin{equation*}
c_{ij} = (H_i+H_j)^{-1}\left(- \lambda_{ij} + \mu_{ij} - S(\alpha^i) - S(\alpha^j) \right)
\end{equation*}

If no constraint is binding, that is if $\alpha$ is an unconstrained solution, $\lambda_{ij}=0 $ and $\mu_{ij}=0 $, and we have the same threshold as in point (1). If we have that $\alpha$ is a constrained solution ($\alpha_{i,j} = 0 $), we have that $\mu_{ij} > 0 $ and as $\forall (i,j) \in (\Iintv*{1, N})^2, H_i+H_j$ is increasing, every value of $c_{i,j}$ that is under the threshold of (1) is possible.
\end{proof}

\begin{prop}\label{app_prop_rebuildP}
For any compatibility list $C$ of dimension $n_I$, there exists $\mathscr{I}$ an individual set of the space $\mathcal{I}^I = (\mathcal{E}, \innerproduct{.}{.}) \times \mathcal{U}_I$ with $\mathcal{E}$ an euclidean space of dimension $K\in\mathbb{N}^{*}$. Moreover, $K \geq I$ is a sufficient condition for $\mathcal{E}$ to exist.
\end{prop}

\begin{proof}[Proof (Proposition \ref{app_prop_rebuildP})]
We first build a compatibility matrix $\tilde{C} \in \mathcal{M}(I,I)$ from the compatibility list $C$ (and denoting $P_i$ the characteristic vector in $\mathcal{E}$ of the individual $i$):
$$\left(\tilde{C}\right)_{i,j} = \begin{cases}
    c_{i,j} = \innerproduct{P_i}{P_j} &\text{if } i \neq j\\
    \eta_i &\text{else}
\end{cases}$$

First, $\tilde{C}$ is symmetric by definiton. Every real symmetric matrix has at least one positive definite submatrix. Let's denote $L$ the family of vectors $P$ such that the submatrix $\tilde{C}_L$ of their scalar product is positive definite. We take $\bar{L} \in \text{sup}\left\{L : \tilde{C}_L \succ 0 \right\}$ and denote $\bar{l} = \text{card}(\bar{L})$. 

Second, we claim that if $\tilde{C}_{\bar L} \succ 0$ ($\tilde{C}_{\bar L}$ is positive definite), then the $\bar L$ family of vectors $P$ is linearly independent. We prove that claim by contradiction. Let's assume that we could have $\tilde{C}_{\bar L} \succ 0$ and $\bar L$ not free. 
\begin{align*}
    \forall \mathbf{x} \in \mathbb{R}^{\bar l}, \mathbf{x}^t\tilde{C}_{\bar L}\mathbf{x} &= \sum_{i=1}^{\bar l}\sum_{j=1}^{\bar l}x_i\left(\tilde{C}_{\bar L}\right)_{i,j}x_j\\
    &= \sum_{i=1}^{\bar l} \eta_i + \sum_{i=1}^{\bar l}\sum_{j=1, i\neq j}^{\bar l}\innerproduct{x_iP_i}{x_jP_j}
\end{align*}

If we define $\eta_i =\innerproduct{P_i}{P_i}$, we have $\mathbf{x}^t\tilde{C}_{\bar L}\mathbf{x} = \left\vert\left\vert  \sum_{i=1}^{\bar l}x_iP_i \right\vert\right\vert^2$. 
$$\tilde{C}_{\bar L} \succ 0 \Leftrightarrow \forall \mathbf{x} \in \mathbb{R}^{\bar l}, \mathbf{x}^t\tilde{C}_{\bar L}\mathbf{x} > 0 \Leftrightarrow \forall \mathbf{x} \in \mathbb{R}^{\bar l}, \left\vert\left\vert  \sum_{i=1}^{\bar l}x_iP_i \right\vert\right\vert^2 > 0$$

But if $\bar L$ is linearly dependant, there exists $\mathbf{x}$ such that $\left\vert\left\vert  \sum_{i=1}^{\bar l}x_iP_i \right\vert\right\vert^2 = 0$ breaking the positive definiteness condition.

Finally, as every vector of $\bar L$ is a vector of $\mathbb{R}^K$, we need $K$ to be at least $\bar l = \text{card}(\bar L)$. By definition, we have $K \geq I \geq \bar l$.
\end{proof}

\begin{defi}[Pairwise stability, from (Jackson 1996)]\label{app_def_pairwiseStable}
A graph $\mathcal{G}(\mathscr{I} , \alpha, C)$ is pairwise stable, iff
\begin{enumerate}
    \item for all $(i,j)$ such that $\alpha_{i,j}$ is a non-null component of $\alpha$, we have $$\forall k \in \{i,j\}, \quad U_k(\alpha, C) \geq U_k(\alpha, C)\vert_{\alpha_{i,j}=0 } $$
    \item for all $(i,j)$ such that $\alpha_{i,j}$ is null component of $\alpha$, we have $$\forall a \in (0 ,1), \forall k \in \{i,j\}, \quad U_k(\alpha, C) \geq U_k(\alpha, C)\vert_{\alpha_{i,j}=a} $$
\end{enumerate}
\end{defi}

\begin{prop}\label{app_prop_pairwiseStable}
Every graph $\mathcal{G}(\mathscr{I} , \alpha, C)$ that is solution of the Social Optimization Problem (Definition \ref{app_def_SOP}) is pair-wise stable. 
\end{prop}

\begin{proof}[Proof (Proposition \ref{app_prop_pairwiseStable})]
We prove both results (1 and 2) by contradiction. Assume that there could exist an $\alpha'$ in the feasible set such that $\exists (i,j) \in \Iintv*{1, N }^2, \forall k \in \{i,j\},\, U_k(\alpha', C) \geq U_k(\alpha, C)$ and $\forall k \in \Iintv*{1, N }^2\backslash\{(i,j), (j,i)\}, \alpha_k = \alpha'_k$. By construction of utility functions, we have that $\forall k \in \Iintv*{1, N }\backslash\{i,j\}, U_k(\alpha, C) = U_k(\alpha', C)$. Hence, we have $\sum_{k\leq N} U_k(\alpha', C) \geq \sum_{k\leq N} U_k(\alpha, C) \Rightarrow \mathcal{O}(\alpha', C) \geq \mathcal{O}(\alpha, C)$. As $\alpha$ is the solution of the SOP, and as the solution is unique, this inequality is absurd. 
\end{proof}

\subsection{Alky utility function}

\begin{defi}[ALKY utility function]\label{app_def_ALKYUtility}
We consider the set $\mathcal{U}^{\text{ALKY}}_I$ of applications $U : [0 ,1)^{I-1} \times (\mathbb{R}^{*}_{+})^{I-1} \longrightarrow \mathbb{R}$:
$$\mathcal{U}^{\text{ALKY}}_I = \left\{\kappa > 1, \gamma  >1, \delta >1, U : \left(\alpha, c \right) \longmapsto \sum^{I-1}_{j=1} \left[ \kappa \alpha_{j} c_{j} - \frac{\alpha_{j}^{\gamma}}{1-\alpha_{j}^{\gamma}} \right] - \left( \sum^{I-1}_{j=1} \alpha_{j} \right)^{\delta}\right\}$$
\end{defi}

\begin{rmk}
Analytically, this utility function is made of a three parts: a positive one which weigths affinity $(c_{i, j})_{i,j \leq N}$ by the relationship intensity $(\alpha_{i, j})_{j \leq N}$ and two costs, that convexly penalize intensity, and create substitution effects between relationships.

$\gamma$ is the convexity parameter of relationships' direct cost (and prevents corner solutions). $\delta$ enable the substitutability effects between relationships. For $\delta > 1$ relations are substitutable: an agent prefers reducing the emphasis on relationships with lower affinity in order to increase the weight on his most successful matches. 
\end{rmk}

\begin{thm}[ALKY Utility]\label{app_thm_utilityALKY}
$\mathcal{U}^{\text{ALKY}}_I$ is a subset of $\mathcal{U}_I$, that is all ALKY functions are utility functions.
\end{thm}

\begin{proof}[Proof (Theorem \ref{app_thm_utilityALKY})]
Let's take $U \in \mathcal{U}^{\text{ALKY}}_I$, $i \in \Iintv*{1, I-1 }$ and $j \in \Iintv*{1, I-1 }, j\neq i$. 

First, (1) and (2) are trivially met $\forall C \in \mathbb{R}^{n_I}, U(0 ,C) = 0 $ and $\nabla U(0 ) > 0 $ iff $C$ is a strictly positive vector.

$$\frac{\partial U}{\partial \alpha_{i}} =  \underbrace{\kappa c_{i} - \frac{\gamma \alpha_{i}^{\gamma-1}}{(1-\alpha_{i}^{\gamma})^2}}_{H(c_i, \alpha_i)} \underbrace{- \delta \left( \sum^{I-1}_{k=1} \alpha_{k} \right)^{\delta -1}}_{S(\alpha)}$$

\begin{align*}
\frac{\partial^2 U}{\partial \alpha_{i} \partial \alpha_{j}} &= \underbrace{-\delta(\delta -1)\left( \sum^{I-1}_{k=1} \alpha_{k} \right)^{\delta -2}}_{s(\alpha)} < 0  \\
\frac{\partial^2 U}{\partial \alpha_{i}^2} &= \underbrace{\frac{\gamma \alpha_{i}^{\gamma-2}}{(1-\alpha_{i}^{\gamma})^3}\left[(1-\gamma) - (1+\gamma)\alpha_{i}^{\gamma} + 2\gamma \alpha_{i}^{2\gamma} \right]}_{h(c_i, \alpha_i)}  \underbrace{- \delta (\delta -1)\left( \sum^{I-1}_{k=1} \alpha_{k} \right)^{\delta -2}}_{s(\alpha)}
\end{align*}

If $s(\alpha)$ is trivially strictly negative for all non-zero vector, we need to consider the sign of $h(c_i, \alpha_i)$. As $\gamma > 1, \delta > 1$, the only remaining term of undetermined sign is: $(1-\gamma) - (1+\gamma)\alpha_{i}^{\gamma} + 2\gamma \alpha_{i}^{2\gamma}$. We switch variables by defining $X = \alpha_{i}^{\gamma}$ and get the following second order polynomial: $(1-\gamma) - (1+\gamma)X + 2\gamma X^2$ evaluated on [0 ,1). It is easy to observe that the polynomial is convex, and is negative on [0 ,1). Hence, $h(c_i, \alpha_i) < 0 $.

\end{proof}

\begin{rmk}
Any individuals set that is only composed of ALKY utility functions satisfies the conditions of Theorem \ref{app_thm_alphaIsAlwaysOpti}: as ALKY functions are linear in C, the sum of their partial derivative with respect to $\alpha$ is trivially bijective.
\end{rmk}

\newpage

\section{Detailed algorithms for the agent-based model}

Algorithms \ref{app_alg:agentTrueAlgo} and \ref{app_alg:agent_responseTrueAlgo} are provided as they are closer to the implementation in C++. $\lambda = 0 .1$ is the learning rate, $\Omega = 0 .1$ the max step size, $\omega=0 .05$ is the minimum weight (above which the link cannot exist) and $\varpi = 0 .0001$ the minimum accepted gradient value. In the C++ implementation, the weights are computed using sigmoid, and the gradient moove is made on the argument of sigmoid and not directly on the weight. 

\begin{figure}[h]
\begin{algorithm}[H]
\caption{Individual $i$ select action protocol (Detailed)}\label{app_alg:agentTrueAlgo}
\begin{algorithmic}
\Function{takeAction}{$i$}
\State $\mathcal{S}^i \gets \{j : d(j, i) \leq \delta \}$ \Comment{Get individuals in the scope of $i$}

\State $C^i_\mathcal{S} \gets \{c_{i,j} : j \in \mathcal{S}^i$ \} \Comment{Get the set of compatibilities}

\State $\alpha^i_\mathcal{S} \gets \{\alpha_{i,j} : j \in \mathcal{S}^i$ \} \Comment{Get the set of weights}

\For{$j \in \mathcal{S}^i $} 
    \If{$j \in M$} \Comment{Check the memory set}
        \State $\nabla U^i_{(j)}(\alpha^i_\mathcal{S}, C^i_\mathcal{S}) \gets 0 $
    \ElsIf{$\nabla U^i_{(j)}(\alpha^i_\mathcal{S}, C^i_\mathcal{S}) \leq \varpi$}
        \State $\nabla U^i_{(j)}(\alpha^i_\mathcal{S}, C^i_\mathcal{S}) \gets 0 $
    \EndIf
\EndFor

\If{$\text{card}(M) > 3$}
    \State \Call{ForgetAMemory}{}  \Comment{Delete one memory}
\EndIf

\State $k \gets \text{softargmax}_j |\nabla U^i(\alpha^i_\mathcal{S}, C^i_\mathcal{S})|$ 

\If{$\nabla U^i_{(k)} > 0 $}
    \State $M \gets M \cup k$ \Comment{Remember the action}
    \State \Return \Call{request}{$k, i, \alpha_{i,k} + \lambda \nabla U^i_{(k)}$} \Comment{Ask $k$ to increase weight}
\ElsIf{$\nabla U^i_{(k)} < 0 $} 
    \State $M \gets M \cup k$ \Comment{Remember the action}
    \State $\alpha_{i,k} \gets \alpha_{i,k} + \lambda \nabla U^i_{(k)}$ \Comment{Direct decrease in friendship with $N^k$}
    \State \Return true
\Else
    \State \Return false \Comment{No action has been taken}
\EndIf
\EndFunction
\end{algorithmic}
\end{algorithm}

\begin{algorithm}[H]
\caption{Individual $k$ response protocol to $i$ (request function, detailed)}\label{app_alg:agent_responseTrueAlgo}
\begin{algorithmic}
\Function{request}{$k, i, a$}
    \State $\mathcal{S}^k \gets \{j : d(j, k) \leq \delta \} \bigcup \{i\}$ \Comment{Get individuals in the scope}
    
    \State $C^k_\mathcal{S} \gets \{c_{k,j} : j \in \mathcal{S}^k$ \} \Comment{Get the set of compatibilities}
    
    \State $\alpha^i_\mathcal{S} \gets \{\alpha_{k,j} : j \in \mathcal{S}^k$ \} \Comment{Get the set of weights}
        
    \If{$\nabla U^k_{(i)}(\alpha^k_\mathcal{S}, C^k_\mathcal{S}) > 0  \wedge \text{min}(a, \alpha_{i,k} + \lambda \nabla U^k_{(i)}) > \omega$}
        \If{$\alpha_{i,k} > 0 $}
            \State $\alpha_{i,k} \gets \text{min}(\text{min}(a, \alpha_{i,k} + \lambda \nabla U^k_{(i)}),\alpha_{i,k} +\Omega)$
        \Else 
            \State $\alpha_{i,k} \gets \Omega$
        \EndIf
        \State \Return true
    \ElsIf{$\nabla U^k_{(i)}(\alpha^k_\mathcal{S}, C^k_\mathcal{S}) < 0 $} \Comment{If he wants to reduce the weight}
        \State $\alpha_{i,k} \gets \text{max}(\text{max}(\alpha_{i,k} + \lambda \nabla U^k_{(i)},\alpha_{i,k} -\Omega), 0 )$ 
        \State \Return true
    \EndIf
    \State \Return false
\EndFunction

\end{algorithmic}
\end{algorithm}
\end{figure}

\begin{figure}[h]
\begin{algorithm}[H]
\caption{Simulation body (detailed)}\label{app_alg:mainTrueAlgo}
\begin{algorithmic}
    \State $t \gets 1$ \Comment{Initialize the time (round counter)}
    \State $v \gets 0 $ \Comment{Initialize an inactive rounds counter}

    \For{$t \leq \text{MAX-SIMULATION-TIME} \wedge v \leq 50+3 $} 
        \State $v_t \gets 0 $ 
        \For{$i \in o_t \in \mathscr{S}(\Iintv*{1, N }) $} \Comment{For each individuals in a random order}
            \If{\Call{takeAction}{i} = false} \Comment{If no action has been taken}
                \State $v_t \gets v_t +1$
            \EndIf
        \EndFor
        \If{$v_t = I$} 
            \State $v \gets v +1$
        \EndIf
        \If{$t > \text{MAX-SIMULATION-TIME} \times 10\%$} 
            \For{$\alpha \in \{\alpha_{i,j} < \omega : i \leq I, j \leq I \}$} 
                \State $\alpha \gets 0 $ \Comment{Removal of the smallest edges}
            \EndFor
        \EndIf
    \EndFor \Comment{End of the simulation process}
\end{algorithmic}
\end{algorithm}
\end{figure}

\end{document}